\documentclass[11pt]{article}
\usepackage{amsmath,amsbsy,amsfonts,amssymb}

\usepackage[frenchb]{babel}
 \newtheorem{thm}{Th\'eor\`eme}
 
 \newtheorem{cor}[thm]{Corollaire}
   
 \newtheorem{lem}[thm]{Lemme}
  
 \newtheorem{prop}[thm]{Proposition}
   
 \newtheorem{rmq}[thm]{Remarque}

   \newcommand{\GL}{\rm GL}
   \newcommand{{\Speh}}{\rm Speh}
   \newcommand{\Jord}{\rm Jord}

         \newcommand{\ord}{\rm ord}
\begin{document}

\title{Fonctions L de paires pour les groupes classiques}
\author{C. M{\oe}glin\\
Institut de Math\'ematiques de Jussieu\\
CNRS, 4 place Jussieu, F-75005 Paris}
\date{}
\maketitle
\abstract{Il y a au moins deux fa\c{c}ons naturelles d'associer des fonctions $L$ de paires \`a un couple form\'e d'une repr\'esentation cuspidales d'un groupe lin\'eaire et d'une repr\'esentation automorphe de carr\'e int\'egrable (irr\'eductible) d'un groupe classique. Les deux d\'efinitions s'appuient sur les travaux tr\`es r\'ecents d'Arthur. La premi\`ere fa\c{c}on est d'utiliser le produit sur toutes les places des facteurs $L$ de paires d\'efinis avec les param\`etres de Langlands. La deuxi\`eme fa\c{c}on de faire, consiste \`a utiliser le transfert endoscopique tordu au groupe lin\'eaire convenable et d'utiliser les fonctions $L$ de paires des groupes lin\'eaires.

Dans cet article on compare ces deux types de fonctions $L$. On montre que celles d\'efinies avec les param\`etres de Langlands ont moins de p\^oles celles d\'efinies par transfert endoscopique. On montre aussi que l'existence de p\^oles pour les fonctions $L$ obtenues avec les param\`etres de Langlands donnent des conditions suffisantes pour l'existence de p\^oles \`a certaines s\'eries d'Eisenstein alors que dans la m\^eme situation  l'existence de p\^oles pour les fonctions $L$ d\'efinies par transfert endoscopiques   sont, elles, des conditions n\'ecessaires pour l'existence de p\^oles \`a ces s\'eries d'Eisenstein.}
\tableofcontents
\section{Notations \label{notations}}
On fixe $k$ un corps de nombres et  $G$ un groupe classique auquel s'applique les travaux d'Arthur, c'est-\`a-dire un groupe special orthogonal ou symplectique. On note $m^*_{G}$ la dimension de la repr\'esentation naturelle du $L$-groupe de $G$.

Soit $d$ un entier et $\rho$ une repr\'esentation cuspidale irr\'eductible de $\GL(d)$; en g\'en\'eral on fixera plut\^ot d'abord $\rho$ et on notera $d$ par $d_{\rho}$ ce qui d\'efinit $d_{\rho}$. On \'ecrit $\vert\, \vert$ pour $\vert det_{\GL(d_{\rho},{\mathbb A})}\vert$ o\`u l'on consid\`ere la valeur absolue ad\'elique du d\'eterminant. Soit $b$ un entier, on d\'efinit ${\Speh}(\rho,b)$ la repr\'esentation de $\GL(ad_{\rho})$ qui se r\'ealise dans l'espace des r\'esidus des s\'eries d'Eisenstein:
$$
\biggl(\prod_{i\in [1,a-1]}(s_{i}-s_{i+1}-1)E(\rho\vert\, \vert^{s_{1}}\times \cdots\times  \rho\vert\, \vert^{s_{a}},f )\biggr)_{s_{1}=(a-1)/2, \cdots, s_{a}=-(a-1)/2},
$$
o\`u $f$ parcourt l'ensemble des sections du fibr\'e localement constant des induites $\rho\vert\, \vert^{s_{1}}\times \cdots\times  \rho\vert\, \vert^{s_{a}}$ o\`u $s_{1}, \cdots, s_{a}\in {\mathbb C}^a$. Cette d\'efinition a une variante locale, o\`u pour toute place $v$ de $k$, ${\Speh}(\rho_{v},b)$ est la composante locale de la repr\'esentation pr\'ec\'edente en la place $v$; on peut en donner une d\'efinition purement locale comme l'unique quotient irr\'eductible de l'induite $\rho_{v}\vert\,\vert^{(a-1)/2} \times \cdots \times \rho_{v}\vert\,\vert^{-(a-1)/2}$.

Pour $\pi$ une repr\'esentation de carr\'e int\'egrable de $G$ et pour $s\in {\mathbb C}$, on consid\'erera aussi les s\'eries d'Eisenstein $E(\rho\times \pi,s)$; elles sont relatives \`a un groupe de m\^eme type que $G$ mais si on \'ecrit $G=Aut(V)$ o\`u $V$ est un espace symplectique ou orthogonal, le groupe en question est le groupe des automorphismes de $V$ auquel on a ajout\'e $d_{\rho}$ plans hyperboliques.

\section{Introduction\label{introduction}}
Soit $\pi$ une repr\'esentation automorphe de carr\'e int\'egrable irr\'educ\-tible de $G$; gr\^ace aux travaux annonc\'es d'Arthur \cite{clay}, on sait associer \`a $\pi$ une repr\'esentation de $\GL(m^*_{G},{\mathbb A})$, not\'ee $\pi^{{\GL}}$ avec la propri\'et\'e que pour   toute place $v$ o\`u toutes les donn\'ees sont non ramifi\'ees, les composantes locales $\pi_{v}$ et $\pi^{{\GL}}_{v}$ se correspondent par la fonctorialit\'e de Langlands associ\'ee \`a l'inclusion naturelle du $L$-groupe de $G$ dans ${\GL}(m^*_{G},{\mathbb C})$.

Soit $\rho$ une repr\'esentation cuspidale unitaire d'un groupe ${\GL}(d_{\rho})$; depuis (\cite{JS}) on sait d\'efinir la fonction $L(\rho\times \pi^{{\GL}},s)$.

Avec   les r\'esultats locaux d'Arthur \cite{clay}, \cite{annonce}, dans le cas des places finies, on conna\^{\i}t aussi la classification de Langlands des repr\'esentations temp\'er\'ees \`a l'aide de l'endoscopie tordue. On peut donc aussi associ\'e \`a toute repr\'esentation irr\'eductibles en la place $v$ son param\`etre de Langlands et en d\'eduire un facteur $L$ local, $L(\rho_{v}\times \pi_{v},s)$. Aux places archim\'ediennes, cela est aussi faisable depuis longtemps. Pour presque toute place $v$, le facteur $L(\rho_{v}\times \pi_{v},s)$ co\"{\i}ncide avec le facteur $L(\rho_{v}\times \pi^{{\GL}}_{v},s)$. On peut donc d\'efinir en tant que fonction m\'eromorphe de $s$:
$$
L(\rho\times \pi,s):=\prod_{v}L(\rho_{v}\times \pi_{v},s)
$$
ce qui vaut $L(\rho\times \pi^{{\GL}},s)\prod_{v\in S}L(\rho_{v}\times \pi_{v},s)/L(\rho_{v}\times \pi^{{\GL}}_{v},s)^{-1},$ o\`u $S$ est un ensemble fini de places suffisamment grand pour contenir les places archim\'ediennes et toutes les places o\`u il y a de la ramification.

Evidemment en g\'en\'eral, $L(\rho\times \pi,s)\neq L(\rho\times \pi^{{\GL}},s)$ et la fonction $L$ de Langlands contient des informations un peu diff\'erentes de la fonction obtenue par l'endoscopie tordue car elle voit les mauvaises places. En particulier, elle a en g\'en\'eral beaucoup moins de p\^oles. Le but de cet article est de rendre pr\'ecise une telle affirmation.

Il y a une difficult\'e \`a comprendre les paquets d'Arthur aux places archim\'ediennes, m\^eme le cas des repr\'esentations temp\'er\'ees n'est pas compl\`etement \'eclairci (cf. \cite{mezo} qui s'appuie sur \cite{shelstad}). Dans tout ce qui suit mais uniquement pour les r\'esultats globaux on suppose que la repr\'esentation $\pi$ a de la cohomologie \`a l'infini et on suppose qu'il existe un nombre r\'eel strictement positif $s_{1}$  tel qu'en toute place $v$ archim\'edienne, l'induite $\rho_{v}\vert\, \vert_{v}^{s_{1}}\times \pi_{v}$ a un caract\`ere infinit\'esimal entier et r\'egulier. 
On montre alors pour tout nombre r\'eel $s_{0}\geq s_{1}$:
$$
\ord_{s=s_{0}}L(\rho\times \pi,s)\geq \ord_{s=s_{0}}L(\rho\times \pi^{\GL},s),
$$
$$
\ord_{s=s_{0}}\bigl(L(\rho\times \pi,s)/L(\rho\times \pi,s+1)\bigr)\geq$$
$$ \ord_{s=s_{0}}\bigl(L(\rho\times \pi^{\GL},s)/L(\rho\times \pi^{\GL},s+1)\bigr).
$$
On commence par d\'emontrer l'analogue de ces r\'esultats d'un point de vue local. Et on montre m\^eme qu'avoir une \'egalit\'e plut\^ot qu'une in\'egalit\'e est \'equivalente \`a une description simple de l'image d'op\'era\-teurs d'entrelacement (cf. \ref{localfini}).

La partie la moins pr\'evisible de cet article est que l'existence d'un p\^ole \`a la fonction $L(\rho\times \pi,s)/L(\rho\times \pi,s+1)$ en $s=s_{0}\geq s_{1}$  est une condition suffisante pour que les s\'eries d'Eisenstein $E(\rho\times \pi,s)$ aient un p\^ole en $s=s_{0}$, au moins si $\pi$ est une repr\'esentation cuspidale. Le r\'esultat n'\'etait pas, de mon point de vue compl\`etement pr\'evisible, car ce qui est clair (cf. \cite{manuscripta}) c'est que l'existence d'un p\^ole en $s=s_{0}$ des s\'eries d'Eisenstein pr\'ec\'edentes entra\^{\i}ne l'existence d'un p\^ole pour la fonction $L(\rho\times \pi^{\GL},s)/L(\rho\times \pi^{\GL},s+1)$ mais il n'y a pas \'equivalence. Pour transformer la condition suffisante que l'on a trouv\'ee en condition n\'ecessaire et suffisante, il faut d\'efinir en toute place $v$, le sous-quotient de Langlands de l'induite $\rho_{v}\vert\,\vert^{s_{0}}\times \pi_{v}$ et on montre alors l'\'equivalence (pour $s_{0}$ avec les hypoth\`eses pr\'ec\'edentes):

\

\noindent
\sl on suppose que $\pi$ est une repr\'esentation cuspidale avec les hypoth\`eses d\'ej\`a faites aux places archim\'ediennes et alors
la fonction $\frac{L(\rho\times \pi,s)}{L(\rho\times \pi,s+1)}$ a un p\^ole en $s=s_{0}$ si et seulement si les deux conditions ci-dessous sont satisfaites:

les s\'eries d'Eisenstein $E(\rho\times \pi,s)$ ont un p\^ole en $s=s_{0}$;

la repr\'esentation r\'esiduelle ainsi d\'efinie, qui est alors n\'ecessaire\-ment irr\'eductible, est produit tensoriel restreint des  sous-quotients de Langlands d\'efinis   ci-dessus.

\

\rm
Quand on sp\'ecialise au cas o\`u $\rho$ est un caract\`ere quadratique, $G$ est un groupe sp\'ecial orthogonal et en supposant que $\pi$ est une repr\'esentation cuspidale, on montre que l'existence d'un p\^ole en $s=s_{0}$ de la fonction $L(\rho\times \pi,s)/L(\rho\times \pi,s+1)$ est une condition suffisante pour que $\pi$ soit dans l'image de s\'eries theta convenables. Ces r\'esultats sont l'objet des paragraphes \ref{eisenstein} et \ref{theta}, bien s\^ur pour les raisons expliqu\'ees ci-dessus on a les hypoth\`eses pr\'ec\'edentes aux places archim\'ediennes.

\

Quand on parle de fonctions $L$, il est in\'evitable de se pr\'eoccuper de l'\'equation fonctionnelle; donc l'article commence par \'etudier les facteurs $\gamma$ le but \'etant de montrer que ce sont les m\^emes que l'on utilise, en la place locale $v$, la composante locale $\pi_{v}$ de $\pi$ ou la composante locale $\pi^{\GL}_{v}$ de $\pi^{\GL}$. Quand $v$ est une place finie, on obtient facilement le r\'esultat car on conna\^{\i}t assez bien les paquets d'Arthur locaux. Quand $v$ est une place archim\'edienne, cela se r\'ev\`ele beaucoup plus difficile et on est oblig\'e de mettre des hypoth\`eses; on renvoit le lecteur \`a \ref{gammaarchimedien} pour ces hypoth\`eses. Quand ces hypoth\`eses sont satisfaites, on obtient alors facilement que $L(\rho\times \pi,s)$ v\'erifie la m\^eme \'equation fonctionnelle que $L(\rho\times \pi^{\GL},s)$, ici pour tout $s\in {\mathbb C}$.

\

Comme expliqu\'e ci-dessus, nos d\'efinitions des fonctions $L$ se fait via le transfert endoscopique de \cite{clay} et \cite{annonce} \`a l'aide des d\'efinitions pour les groupes $\GL$. Gr\^ace aux r\'esultats de Shahidi \cite{shahidifini}, \cite{shahidiarchimedien} et d'Henniart \cite{henniart} cela revient, quand on travaille localement, \`a utiliser les d\'efinitions du c\^ot\'e galoisien. Globalement, on se ram\`ene ainsi aux travaux de Jacquet et Shalika sur la d\'efinition et l'\'etude des fonctions $L$ de paires \cite{JS}.

Quand $\rho$ (ci-dessus) est un caract\`ere les fonctions $L$ et les facteurs $\epsilon$ ont \'et\'e d\'efinis, uniquement en termes des repr\'esentations par Rallis et Piatetskii-Shapiro en \cite{RPS}, chapitre 1; localement les d\'efinitions sont donn\'ees par Lapid et Rallis en \cite{LR}. Comme Rallis et Soudry en \cite{RS} ont montr\'e les propri\'et\'es de stabilit\'e de ces facteurs locaux, les m\'ethodes d'Henniart s'appliquent pour v\'erifier que les facteurs locaux d\'efinis par ces travaux sont les m\^emes que ceux utilis\'es ici, au moins aux places finies,  gr\^ace \`a l'\'equation fonctionnelle d\'emontr\'ee en \ref{gammafini}. Aux places archim\'ediennes, le r\'esultat reste vrai sous les hypoth\`eses de \ref{gammaarchimedien}. Du point de vue global,  il reste encore du travail \`a faire pour avoir vraiment la co\"{\i}ncidence.

\

Pour finir l'article, on s'int\'eresse au nombre de p\^oles des fonctions $L(\rho\times \pi,s)/L(\rho\times \pi,s+1)$ quand $s$ parcourt ${\mathbb R}_{>0}$; pour avoir un r\'esultat satisfaisant, on suppose qu'en toute place archim\'edienne le groupe est quasi d\'eploy\'e et que la composante locale de $\pi$ est dans le paquet de Langlands associ\'e au paquet d'Arthur. Avec cette hypoth\`ese, en les places archim\'ediennes, les deux d\'efinitions de facteurs $L$ co\"{\i}ncident et tous les r\'esultats pr\'ec\'edents sont donc vrais sans hypoth\`ese supp\'ementaire. On suppose encore que $\pi$ est une repr\'esentation cuspidale et on \'ecrit $\pi^{\GL}$ comme une induite de repr\'esentations de carr\'e int\'egrable, $\pi^{\GL}\simeq \times_{(\rho',b')\in {\mathcal {S}}}{\Speh}(\rho',b')$ (ce qui d\'efinit ${\mathcal {S}}$). Alors on montre que le nombre de p\^oles cherch\'e est inf\'erieur ou \'egal au cardinal des couples $(\rho',b')\in {\mathcal {S}}$ tel que $\rho'\simeq \rho$ et $(\rho',b'+2)\notin {\mathcal{S}}$. Sans l'hypoth\`ese que $\pi$ est cuspidal le nombre de p\^ole de la fonction consid\'er\'e est inf\'erieur ou \'egal \`a deux fois le cardinal des couples $(\rho',b')\in {\mathcal {S}}$ tel que $\rho'\simeq \rho$ et on s'attend \`a ce que cette borne puisse \^etre atteinte.

\

Je remercie Nicolas Bergeron pour avoir attir\'e mon attention sur la caract\'erisation de l'image des s\'eries theta et Jean-Loup Waldspurger pour les discussions que nous avons eues.

\section{Propri\'et\'es des facteurs $\gamma$\label{gamma}}
On fixe $\pi$ une repr\'esentation de carr\'e int\'egrable de $G$ et on note $\pi^{\GL}$ la repr\'esentation automorphe de $\GL(m^*_{G},k)$ qu'Arthur lui associe. On fixe aussi $\rho$ une repr\'esentation cuspidale irr\'eductible et unitaire de $\GL(d_{\rho},k)$ (ce qui d\'efinit $d_{\rho}$)
\subsection{Le cas des places finies\label{gammafini}}
On fixe une  place finie $v$ de $k$. On consid\`ere $\pi_{v}$ comme quotient de Langlands d'une induite:
$$
 \times_{i\in [1,\ell]}\sigma_{i}\vert\,\vert^{-s_{i}}\times \pi_{temp}\rightarrow \pi_{v},
$$
o\`u pour tout $i\in [1,\ell]$, $\sigma_{i}$ est une s\'erie discr\`ete unitaire irr\'eductible, $s_{i}\in {\mathbb R}_{>0}$, $\pi_{temp}$ est une repr\'esentation temp\'er\'ee et $s_{1}\geq \cdots \geq s_{\ell}$. On note $\pi^{{\GL}}_{temp}$ la repr\'esentation du groupe lin\'eaire convenable qui correspond \`a $\pi_{temp}$ par l'endoscopie tordue. 

Soit $\rho_{v}$ la composante locale de $\rho$ et soit $s\in {\mathbb C}$. On sait d\'efinir, \`a la suite de Shahidi (\cite{shahidifini}) les facteurs $\gamma$ dans les groupes lin\'eaire et on pose:
$$
\gamma(\rho_{v},\pi_{v},s):=\times_{i\in [1,\ell]}\gamma(\rho_{v}\times \sigma_{i}^*,s+s_{i})\gamma(\rho_{v}\times \sigma_{i}^*,s-s_{i}) \times \gamma(\rho_{v}\times \pi_{temp}^{{\GL}},s).
$$
Evidemment avec des d\'efinitions analogues, on d\'efinit les facteurs $L$ et les facteurs  $\epsilon(\rho_{v}\times \pi_{v},s)$ qui intervient dans l'\'equation fonctionnelle, ci-dessous.

En utilisant $\pi^{\GL}_{v}$ au lieu de $\pi_{v}$, on sait aussi d\'efinir \`a la suite de Shahidi
$
\gamma(\rho_{v} \times \pi(\psi_{v})^{{\GL}},s)$.

\begin{prop} On a l'\'egalit\'e de fonction m\'eromorphe: $\gamma(\rho_{v},\pi_{v},s)=\gamma(\rho_{v}\times \pi(\psi_{v})^{{\GL}},s)$.
\end{prop}
C'est un probl\`eme purement dans les groupes lin\'eaires.

Les facteurs $\gamma$ jouissent de tr\`es bonnes propri\'et\'es r\'esum\'ees dans les dix commandements de \cite{LR}. On peut  d\'efinir $\gamma(\tau_{v}\times \sigma,s)$, o\`u $\sigma$ est une induite non n\'ecessairement irr\'eductible et  $\gamma(\tau_{v}\times \sigma,s)=\gamma(\tau_{v}\times \sigma',s)$ pour tout sous-quotient irr\'eductible $\sigma'$ de $\sigma$. On pose 
$\sigma=\times_{i\in [1,\ell]}\sigma_{i}\vert\,\vert^{-s_{i}}\times \pi_{temp}\times \sigma_{i}^*\vert\,\vert^{s_{i}}$ et avec la multiplicativit\'e des facteurs $\gamma$, on a:
$$
\gamma(\rho_{v},\pi_{v},s)=\gamma(\rho_{v}\times \sigma,s).
$$
Pour d\'emontrer la proposition, il suffit de v\'erifier que $\pi(\psi_{v})^{{\GL}}$ et $\sigma$ ont m\^eme support cuspidal. Le support cuspidal de $\sigma$ est ce que l'on a appel\'e le support cuspidal \'etendu de $\pi_{v}$ en \cite{surlespaquets} 4.1 et on a montr\'e que que c'est aussi le support cuspidal de $\pi(\psi_{v})^{{\GL}}$ pour tout paquet contenant $\pi$ en \cite{surlespaquets} 4.1 et 4.2. Cela termine la preuve.
\subsection{Le cas des places archim\'ediennes\label{gammaarchimedien}}
Ici on fixe $v$ une place archim\'edienne; si $v$ est complexe, on ne fait aucune hypoth\`ese mais si $v$ est une place r\'eelle, on suppose que $\pi_{v}$ a de la cohomologie et que $\pi^{\GL}_{v}$ param\'etrise un paquet d'Adams-Johnson contenant $\pi_{v}$.  Et la proposition ci-dessus est aussi exacte:

\begin{prop} On a l'\'egalit\'e de fonction m\'eromorphe: $\gamma(\rho_{v},\pi_{v},s)=\gamma(\rho_{v}\times \pi_{v}^{{\GL}},s)$.
\end{prop}
On commence  par le cas o\`u $v$ est  une place r\'eelle. Sans aucunes hypoth\`eses on sait que $\pi^{\GL}_{v}$ est une induite irr\'eductible de repr\'esentations ${\Speh}(\delta,b)$, \'eventuellement tordu par un caract\`ere $\vert det\vert^x$ avec $x$ un nombre r\'eel de valeur absolue strictement inf\'erieure \`a $1/2$ et o\`u $\delta$ est soit une s\'erie discr\`ete de $\GL(2,{\mathbb R})$ soit un caract\`ere unitaire de ${\mathbb R}^*$. On note $\Jord(\pi^{\GL}_{v})$ l'ensemble des triplets $(\delta,b,x)$ qui interviennent ci-dessus.

D\`es que l'on suppose que le caract\`ere infinit\'esimal de $\pi_{v}$ est r\'egulier, il y a beaucoup de simplification et on a encore plus de simplification si l'on suppose que $\pi_{v}$ a un caract\`ere infinit\'esimal entier et r\'egulier. Sous cette derni\`ere hypoth\`ese pour tout $(\delta,b,x)\in \Jord(\pi^{\GL}_{v})$, on a n\'ecessairement $x=0$, $\delta$ est autodual et la parit\'e de $b$ est d\'etermin\'e par $\delta$. Ici on utilise simplement le fait que le caract\`ere infinit\'esimal est entier; on r\'eduit donc $Jord(\pi^{\GL}_{v})$ a un ensemble de couples $(\delta,b)$. La r\'egularit\'e du caract\`ere infinit\'esimal assure que pour $\delta$ une s\'erie discr\`ete de $\GL(2,{\mathbb R})$ fix\'ee, il existe au plus un entier $b$ tel que $(\delta,b)\in \Jord(\pi^{\GL}_{v})$. De plus except\'e pour les groupes orthogonaux pairs, il existe au plus un caract\`ere quadratique $\eta$ tel qu'il existe un entier $b$ avec $(\eta,b)\in \Jord(\pi^{\GL}_{v})$. Si $G(k_{v})$ est un groupe sp\'ecial orthogonal pair, alors s'il existe $(\eta,b)\in \Jord(\pi^{\GL}_{v})$ avec $\eta$ un caract\`ere quadratique n\'ecessairement $b$ est impair et $\Jord(\pi^{\GL},v)$ contient n\'ecessairement exactement  deux tel couples l'un est $(\eta,b)$ avec $b\geq 1$ mais n\'ecessairement impair et l'autre est $(\eta\eta_{G},1)$ avec $\eta_{G}$ est le caract\`ere trivial si $G(k_{v})$ est une forme int\'erieur d'un groupe d\'eploy\'e et est le caract\`ere non trivial sinon.

Par d\'efinition $\pi_{v}$ est un quotient de Langlands d'une repr\'esentation induite de la forme 
$
ind \otimes _{i\in [1,\ell]}\delta_{i}\vert det\vert^{x_{i}}\otimes \tau,
$
o\`u les $\delta_{i}$ sont soit des s\'eries discr\`etes de $\GL(2,{\mathbb R})$ soit des caract\`eres unitaires de ${\mathbb R}^*$ est les $x_{i}$ sont des nombres r\'eels d\'ecroissants strictement positifs et $\tau$ est une repr\'esentation temp\'er\'ee d'un groupe de m\^eme type que $G(k_{v})$. 
On note $\tau^{\GL}$ la repr\'esentation du groupe lin\'eaire convenable qui correspond au $L$-paquet contenant $\tau$ (maintenant on a une bonne r\'ef\'erence avec les travaux de \cite{mezo}). Et l'\'egalit\'e des facteurs $\gamma$ que l'on cherche r\'esultera de la propri\'et\'e suivante: la repr\'esentation $\pi^{\GL}_{v}$ et  le quotient de Langlands de l'induite 
$
ind  \otimes _{i\in [1,\ell]}\delta_{i}\vert det\vert^{x_{i}}\otimes \tau^{\GL}\otimes _{i\in [\ell,1]}\delta_{i}\vert det \vert^{-x_{i}}
$
sont sous-quotient d'une m\^eme s\'erie principale.

L'id\'ee sous-jacente, quand $\pi_{v}$ a de la cohomologie,  est qu'\`a la repr\'esentation $E$ qui est le syst\`eme de coefficients dans lequel $\pi_{v}$ a de la cohomologie, on associe une repr\'esentation irr\'eductible de $GL(m^*_{G},{\mathbb R})$, $E^{\GL}$ et la  s\'erie principale que l'on cherche est une s\'erie principale ayant $E^{\GL}$ comme quotient de Langlands. Faute de r\'ef\'erences et de comp\'etences, on fait un peu diff\'eremment et de fa\c{c}on plus terre \`a terre de plus on utilise ci-dessous l'hypoth\`ese forte que $\pi^{\GL}_{v}$ param\'etrise un paquet d'Adams-Johnson (cf. \cite{aj}) contenant $\pi_{v}$.

En effet cette hypoth\`ese permet de relier les $\delta_{i},x_{i}$ qui interviennent ci-dessus \`a $\Jord(\pi^{\GL}_{v})$; le r\'esultat est le suivant pour tout $(\delta,b)\in Jord(\pi^{\GL}_{v})$ il existe $b'\leq b$ un entier de m\^eme parit\'e que $b$, \'eventuellement $b'=0$ si $b$ est pair tel que l'ensemble des $(\delta_{i},x_{i})$ pour $i\in [1,\ell]$ comme ci-dessus soit \`a l'ordre pr\`es exactement l'ensemble des $$
\cup_{(\delta,b)\in \Jord(\pi^{\GL}_{v})}\cup_{x\in [(b-1)/2,(b'+1)/2}(\delta,x),
$$
o\`u si $b'=b$, l'ensemble attach\'e \`a $(\delta,b)$ est vide. De plus $\tau$ est dans le paquet d'Adams-Johnson associ\'e \`a la repr\'esentation induite des ${\Speh}(\delta,b')$ pour $(\delta,b)$ parcourant $\Jord(\pi^{\GL}_{v})$ (ces assertions r\'esultent des descriptions explicites des repr\'esentations dans un paquet d'Adams-Johnson (cf. \cite{aj} \&3) ainsi que de la description explicite des param\`etres de Langlands des repr\'esenta\-tions ayant de la cohomologie donn\'ee dans \cite{vz} 6.16.

Pour simplifier un peu les notations, on commence par remarquer pour avoir l'assertion pour $\pi_{v}$, il suffit maintenant de l'avoir pour $\tau$. Cela permet de supposer d\`es le d\'epart que $\pi_{v}$ est temp\'er\'ee et on note $\pi^{\GL}_{temp}$ la repr\'esentation temp\'er\'ee de $GL(m^*_{G},{\mathbb R})$ qui param\'etrise le paquet de repr\'esen\-tations temp\'er\'ees contenant $\pi_{v}$ (\cite{mezo}). Ainsi $\pi_{v}$ est \`a la fois dans le paquet d\'etermin\'e par $\pi^{\GL}_{v}$ et dans le paquet d\'etermin\'e par  $\pi^{\GL}_{temp}$. Le point est de d\'emontrer que $\pi^{\GL}_{temp}$ et $\pi^{\GL}_{v}$ sont des sous-quotients d'une m\^eme s\'erie principale. 

 On commence par se ramener au cas o\`u $\tau$ est une s\'erie discr\`ete. S'il n'en est pas ainsi, $G$ est un groupe orthgonal pair et en la place $v$ il est forme int\'erieure de la forme d\'eploy\'ee. Et $\tau=\pi_{v}$ est un sous-module irr\'eductible d'une induite \`a partir d'un caract\`ere quadratique $\eta'$ de ${\mathbb R}^*$ avec une s\'erie discr\`ete $\tau'$. On v\'erifie sur la d\'efinition de $\pi_{v}$ que n\'ecessairement $\eta'$ est le caract\`ere $\eta$ intervenant dans la d\'efinition de $\pi^{\GL}_{v}$ et que $Jord(\pi_{v})$ contient deux copies de $(\eta,1)$; en rempla\c{c}ant $\tau$ par $\tau'$ on supprime ces 2 copies. On suppose donc que $\pi_{v}$ est une s\'erie discr\`ete.

On note $\pi^{\GL}_{temp}$ le param\`etre de Langlands de l'unique paquet de repr\'esentations temp\'er\'ees (ici des s\'eries discr\`etes) contenant $\pi_{v}$.
On proc\`ede comme on l'a fait pour \'ecrire $\pi^{\GL}_{v}$ et on trouve que $\pi^{\GL}_{temp}$ est une induite de s\'erie discr\`ete de $\GL(2,{\mathbb R})$ si $G$ est un groupe symplectique, orthogonal impair ou orthogonal pair mais forme int\'erieure de la forme d\'eploy\'ee. Dans le cas restant il s'ajoute la repr\'esentation de $GL(2,{\mathbb R})$ induite du caract\`ere du tore d\'eploy\'e produit du caract\`ere trivial et du caract\`ere signe. Dans tous les cas  $\pi^{\GL}_{temp}$ est sous-quotient d'une s\'erie principale induite \`a partir du tore diagonal et d'un essentiellement produit de valeur absolue \`a une puissance uniquement d\'etermin\'ee par le caract\`ere infinit\'esimal (le essentiellement veut dire que ceci est vrai sauf dans le cas o\`u $G(k_{v})$ est un groupe orthogonal pair dans la classe du groupe quasi-d\'eploy\'e non d\'eploy\'e ou l'un  des caract\`eres est le signe). La m\^eme assertion est vrai pour $\pi^{\GL}_{v}$ si $\eta=1$. Cela r\`egle donc le cas des groupes symplectiques. Dans le cas des groupes orthogonaux, on se ram\`ene \`a ce cas en tensorisant $\pi_{v}$ par $\eta$ compos\'e avec la norme spinorielle et en tensorisant $\rho_{v}$ par le caract\`ere signe du d\'eterminant. Ainsi on n'a pas chang\'e les facteurs $\gamma$ qui nous int\'eressent et on n'a pas chang\'e $\pi^{\GL}_{temp}$ qui est invariant par cette tensorisation; \'evidemment qu'il n'est pas utile de faire cette tensorisation mais il aurait fallu \'ecrire un peu diff\'eremment la s\'erie principale choisie. Cela termine la preuve pour le cas des places r\'eelles.

\

On suppose maintenant que $v$ est une place complexe; ici il n'y a aucun probl\`eme et pas d'hypoth\`ese \`a faire. On \'ecrit d'abord la repr\'esentation de $GL(m^*,{\mathbb C})$ d\'etermin\'ee par le paquet d'Arthur; c'est n\'ecessairement une induite de caract\`ere preque unitaire (c'est-\`a-dire de partie r\'eelle comprise strictement entre $-1/2$ et $1/2$). Ainsi le facteur $\gamma$ associ\'e est aussi le facteur $\gamma$ associ\'e \`a une s\'erie principale de $G(k_{v})=G({\mathbb C})$ induite \`a partir d'un caract\`ere d'un sous-groupe de Borel uniquement d\'etermin\'e par le caract\`ere infinit\'esimal de $\pi_{v}$. Comme le facteur $\gamma$ associ\'e \`a $\pi_{v}$ via ses param\`etres de Langlands a \'evidemment la m\^eme propri\'et\'e, on obtient l'\'egalit\'e des facteurs $\gamma$ cherch\'ee.

\subsection{Application \`a l'\'equation fonctionnelle\label{fonctionnelle}}
Cette section n'a \'evidemment rien d'original. On fixe $\pi$ et $\rho$ comme dans l'introduction. A toute place $v$ de $k$ qui est r\'eelle, on suppose que $\pi_{v}$ a de la cohomologie dans un bon syst\`eme de coefficients et que le paquet d'Arthur d\'etermin\'e en cette place par $\pi$ co\"{\i}ncide avec un des paquets d'Adams-Johnson contenant $\pi_{v}$.

On rappelle l'\'equation fonctionnelle dans les groupes lin\'eaire: soit $\rho,\rho'$ des repr\'esentations cuspidales d'un groupe lin\'eaire et $s\in{\mathbb C}$ alors il existe une fonction holomorphe inversible $\epsilon(\rho\times \rho',s)$ tel que l'on ait l'\'egalit\'e:
$$
L(\rho\times \rho',s)=\epsilon(\rho\times \rho',s)L(\rho^*\times \rho^{'*},1-s).
$$
On a par d\'efinition
$$
L(\rho\times \pi^{{\GL}},s)=\times_{(\rho',b')\in {\Jord}(\pi)}\times_{k\in[(b'-1)/2,-(b'-1)/2])}L(\rho\times \rho',s+k).$$Le point ici est que pour tout $(\rho',b')$, $\rho'$ est autoduale et le segment $[(b'-1)/2,-(b'-1)/2]$ est centr\'e en $0$ (pour ce qui suit il suffit que $\pi^{{\GL}}$ soit autoduale) et il existe donc une fonction holomorphe inversible de $s$ tel que l'on ait l'\'egalit\'e
$$
L(\rho\times \pi^{{\GL}},s)=\epsilon(\rho\times \pi^{{\GL}},s)L(\rho^*\times \pi^{{\GL}},1-s).
$$
\begin{prop} Avec les hypoth\`eses faites au d\'ebut de ce paragraphe, on a l'\'egalit\'e de fonctions m\'eromorphes:
$$
L(\rho\times \pi,s)=\epsilon(\rho\times \pi,s)L(\rho^*\times \pi,1-s).
$$
\end{prop}
En effet d'apr\`es l'\'egalit\'e des facteurs $\gamma$ locaux montr\'es, on a: 
$$
\epsilon(\rho\times \pi,s)L(\rho^*\times \pi,1-s)/L(\rho\times \pi,s)=$$
$$\epsilon(\rho\times \pi^{{\GL}})L(\rho^*\times \pi^{{\GL}},s)/L(\rho\times \pi^{{\GL}}, s)
=1.$$
\section{P\^oles des facteurs et des fonctions $L$}
\subsection{P\^oles des facteurs $L$ locaux aux places finies\label{localfini}}
On fixe $\pi$ et $\rho$ comme pr\'ec\'edemment et une place $v$ de $k$; on suppose que $v$ est une place finie. On a d\'efini $L(\rho_{v}\times \pi_{v},s)$ et $L(\rho_{v}\times \pi^{\GL}_{v},s)$. On commence par \'etudier la fonction $$
\phi_{v}(\rho_{v},\pi_{v},s):=\frac{L_{v}(\rho_{v}\times \pi_{v},s)/L_{v}(\rho_{v}\times \pi_{v},s+1)
}{L_{v}(\rho_{v}\times \pi^{\GL}_{v},s)/L_{v}(\rho_{v}\times \pi^{\GL}_{v},s+1)}.
$$
Cette fonction est li\'ee \`a l'op\'erateur d'entrelacement standard, entre les induites:
$$M_{v}(\rho_{v},\pi_{v},s):
\rho_{v}\vert\, \vert^{s}\times \pi_{v} \rightarrow \rho_{v}^*\vert\, \vert^{-s}\times \pi_{v}.
$$
On va l'\'etudier pour $s=s_{0}$ r\'eel et sup\'erieur ou \'egal \`a 1/2. Le fait de consid\'erer $s$ r\'eel est une simplification d'\'ecriture puisque l'on peut tordre $\rho_{v}$ par un caract\`ere unitaire. Par contre le fait que $s\geq 1/2$ est une restriction d\^ue au fait que l'on ne conna\^{\i}t pas la conjecture de Ramanujan; c'est une vraie restriction mais qui dispara\^{\i}t pour les r\'esultats globaux que l'on a en vue.

On pose $$N^{\GL}(\rho_{v},\pi_{v},s):= \bigl(L(\rho_{v}\times \pi^{\GL}_{v},s)/L(\rho_{v}\times \pi^{\GL}_{v},s+1)
\bigr)^{-1}M_{v}(\rho_{v},\pi_{v},s).
$$
C'est un op\'erateur d'entrelacement et on a montr\'e en \cite{holomorphie} 3.2 et \cite{image} 5.2  que cet op\'erateur est holomorphe en $s=s_{0}\in {\mathbb R}_{\geq 1/2}$.

On pose aussi
$$
N^{\GL}(\rho_{v},\pi_{v},s):= \bigl(L(\rho_{v}\times \pi_{v})/L(\rho_{v}\times \pi_{v},s+1)
\bigr)^{-1}M_{v}(\rho_{v},\pi_{v},s).
$$
Cet op\'erateur est l'op\'erateur normalis\'e \`a la Langlands-Shahidi, il a de moins bonnes propri\'et\'es d'holomorphie que l'op\'erateur pr\'ec\'edent mais, en un point o\`u il est holomorphe, il est plus facile d'avoir une description de son image.

Pour cela on fixe $s_{0}\in {\mathbb R}_{\geq 1/2}$ et on d\'efinit le sous-quotient de Langlands de l'induite $\rho_{v}\vert\,\vert^{s_{0}}\times \pi_{v}$: $\rho_{v}$ est une induite irr\'eductible de repr\'esentations de Steinberg \'eventuellement tordu par un caract\`ere non unitaire mais dont l'exposant est strictement inclus dans $]-1/2,1/2[$. On note $\tau_{-k}, \cdots \tau_{-1}$ ces repr\'esentations ordonn\'ees de telle sorte que l'exposant aille en d\'ecroissant. D'autre part on \'ecrit $\pi_{v}$ comme quotient de Langlands, ainsi il existe des repr\'esentations de Steinberg $\tau_{1}, \cdots, \tau_{\ell}$ tordu par des exposants non unitaires mais strictement positifs rang\'ees dans l'ordre d\'ecroissant des exposants et une repr\'esentation temp\'er\'ee, $\pi_{temp}$. d'un groupe de m\^eme type que $G(k_{v})$ tel que $\pi_{v}$ soit l'unique quotient irr\'eductible de l'induite du produit tensoriel des $(\tau_{i})_{i\in [1,\ell]}$ avec $\tau$. On  note $\sigma$ une permutation de l'ensemble $[-k,\ell]-\{0\}$ croissante sur $[-k,-1]$ et sur $[1,\ell]$ qui r\'earrange la collection de repr\'esentations $\tau_{-k}\vert\,\vert^{s_{0}}, \cdots, \tau_{-1}\vert\,\vert^{s_{0}},\tau_{1}, \cdots, \tau_{\ell}$ de fa\c{c}on \`a ce que les exposants soient dans l'ordre d\'ecroissant. Evidemment $\sigma$ d\'epend de $s_{0}$ et n'est pas uniquement d\'etermin\'ee. Cela permet de construire une induite en ayant r\'eordonn\'e ces repr\'esentations et en induisant avec la repr\'esentation temp\'er\'ee $\pi_{temp}$ et cette induite a un unique quotient irr\'eductible c'est ce que l'on appelle le sous-quotient de Langlands de $\rho_{v}\vert\,\vert^{s_{0}}\times \pi_{v}$. Pour la suite on appelle $\tau(s)$ la repr\'esentation du groupe $\GL$ convenable induites $\times_{j\in [-k,-1]\cup [1,\ell]}\tau_{\sigma^{-1}(j)}$.

Dans la proposition ci-dessous on appelle ordre d'un op\'erateur m\'eromor\-phe d\'ependant de la variable $s$ en $s=s_{0}$ le plus entier $x$ tel que multipli\'e par $(s-s_{0})^{-x}$ il soit holomorphe en $s=s_{0}$ et non identiquement nul en ce point.
\begin{prop}Soit $s_{0}\in {\mathbb R}_{\geq 1/2}$, alors:

(i)$
ord_{s=s_{0}}\phi(\rho_{v},\pi_{v},s)\geq ord_{s=s_{0}}N^{\GL}(\rho_{v},\pi_{v},s) \geq 0.
$

(ii) $
ord_{s=s_{0}}\phi(\rho_{v},\pi_{v},s)=0$ si et seulement si  $N^{\GL}(\rho_{v},\pi_{v},s_{0})\not\equiv 0$ et a pour image  le sous-quotient de Langlands de l'induite $\rho_{v}\vert\,\vert^{s_{0}}\times \pi_{v}$.
\end{prop}
On montre d'abord que l'ordre de l'op\'erateur en $s=s_{0}$, $N^{L}(\rho_{v},\pi_{v},s)$ est n\'egatif ou nul et que si cet ordre est z\'ero, c'est-\`a-dire que l'op\'erateur est holomorphe, alors son image en $s=s_{0}$ contient, comme sous-quotient, le sous-quotient de Langlands de l'induite $\rho_{v}\vert\,\vert^{s_{0}}\times \pi_{v}$. 
On consid\`ere le diagramme d'op\'erateur d'entrelacement standard (et d'inclusions \'evidentes):
$$
\begin{matrix}
 && &\tau(s)\times \pi_{temp}\\ \\
 &&&{\downarrow}\\ \\
&\rho_{v}\vert\,\vert^{s}\times \pi_{v} &\underset{\iota}\hookrightarrow &\times_{j\in [-k,-1]}\tau_{j}\vert\,\vert^{s}\times_{i\in [1,\ell]}\tau_{i}^*\\
&&&\times \pi_{temp}\\
&\downarrow  &&
\\\\
&\rho_{v}^{*-s}\times \pi_{v}&&\downarrow\\
&\iota\downarrow &&\\\\
&\times_{j\in [-1,-k]}\tau_{j}^{*}\vert\,\vert^{-s}\times_{i\in [1,\ell]}\sigma^*_{i}&\rightarrow &\tau(s)^*\times \pi_{temp}\\
&\times \pi_{temp}&&\\
\end{matrix}
$$
Les deux fl\`eches not\'ees $\iota$ sont \'evidemment semblables, c'est l'inclusion de $\pi_{v}$ dans $\times_{\in [1,\ell]}\tau_{i}^*\times \pi_{temp}$.

Le diagramme n'est commutatif qu'\`a une fonction m\'eromorphe pr\`es et de plus la fl\`eche verticale en bas \`a droite peut ne pas \^etre holomorphe en $s=s_{0}$ car on est amen\'e \`a faire l'\'echange des deux facteurs de l'induite $\tau_{i}^* \times \tau_{j}^{*}\vert\,\vert^{-s}$ chaque fois que $-s_{i}\geq -x_{j}-s_{0}$ o\`u l'on a not\'e $s_{i}$ l'exposant de $\tau_{i}$ et $x_{j}$ celui de $\tau_{j}$. Or  quand on a \'egalit\'e, l'op\'erateur a un p\^ole. Pour \'eviter cela (et pour avoir la commutativit\'e du diagramme) on normalise toutes les fl\`eches \`a la Langlands-Shahidi:  la fl\`eche verticale en bas \`a droite est remplac\'ee en multiplicant l'op\'erateur d'entrelacement standard par l'inverse de la fonction m\'eromorphe:
$$
\prod _{i\in [1,\ell]}(L(\rho_{v}\times \tau_{i},s)/L(\rho_{v}\times \tau_{i},s+1))
$$
$$
\times (L(\rho_{v}\times (\pi_{temp})^{\GL},s)/L(\rho_{v}\times (\pi_{temp})^{\GL},s+1) )$$
$$\prod_{i\in [1,\ell],j\in [-k,-1]; \sigma(j)<\sigma(i)}L(\tau_{i}^*\times \tau_{j}^*,s)/L(\tau_{i}^*\times \tau_{j}^*,s+1).
$$
 Dans ces conditions, le compos\'e des fl\`eches verticales \`a droite est holomorphe en $s=s_{0}$ et a pour image exactement le sous-quotient de Langlands de l'induite $\tau(s)\times \pi_{temp}$.

On normalise aussi la premi\`ere fl\`eche verticale en rempla\c{c}ant l'op\'e\-rateur d'entrelacement standard par  l'op\'erateur $N^L(\rho_{v},\pi_{v},s)$ et la derni\`ere fl\`eche horizontale (qui est purement dans un groupe lin\'eaire) \`a la Langlands-Shahidi; cette fl\`eche est alors not\'e $N(\sigma,s)$ et c'est un op\'erateur holomorphe en $s=s_{0}$.

 Avec ces normalisations, le diagramme est commutatif \`a une fonction holomorphe inversible pr\`es en $s=s_{0}$, fonction qui vient des facteurs $\epsilon$ que l'on a n\'eglig\'es et des facteurs de normalisation  \`a la Langlands-Shahidi de la fl\`eche tout en haut \`a droite que l'on a n\'eglig\'e car ils sont holomorphes inversibles (ceux de la fl\`eche tout en haut \`a droite).

On remarque encore que cette fl\`eche tout en haut \`a droite a son image qui est certainement incluse dans $\rho_{v}\vert\,\vert^{s}\times \pi_{v}$. Ainsi $N(\sigma,s)\circ \iota \circ N^L(\rho_{v},\pi_{v},s)$ est holomorphe en $s=s_{0}$ et calcul\'e en ce point a son image qui contient le sous-module de Langlands de $\tau(s_{0})\times \pi_{temp}$. Si $N^L(\rho_{v},\pi_{v},s)$ est holomorphe en $s=s_{0}$, l'holomorphie de $N(\sigma,s)$ en $s=s_{0}$ force l'op\'erateur $N^L(\rho_{v},\pi_{v},s_{0})$ d'avoir une image non nulle et contenant comme sous-quotient le sous-quotient de Langlands de l'induite $\rho_{v}\vert\,\vert^{s_{0}}\times \pi_{v}$.

On d\'emontre maintenant (i): par d\'efinition, on a l'\'egalit\'e d'op\'era\-teurs m\'eromorphes
$$
\phi(\rho_{v},\pi_{v},s)^{-1}N^{\GL}(\rho_{v},\pi_{v},s)=N^L(\rho_{v},\pi_{v},s).$$
Et en prenant les ordres en $s=s_{0}$:
$$
-ord_{s=s_{0}}\phi(\rho_{v},\pi_{v},s)+ord_{s=s_{0}}N^{\GL}(\rho_{v},\pi_{v},s)=$$
$$ord_{s=s_{0}}N^L(\rho_{v},\pi_{v},s) \leq 0.
$$
D'o\`u (i) puisque l'holomorphie de $N^{\GL}(\rho_{v},\pi_{v},s)$ en $s=s_{0}$ a \'et\'e montr\'ee en \cite{holomorphie} 3.2 et \cite{image} 5.2.

Montrons (ii); on suppose d'abord que $\phi(\rho_{v},\pi_{v},s_{0})\neq 0$. D'apr\`es (i), cela force l'op\'erateur $N^{\GL}(\rho_{v},\pi_{v},s_{0})$ a \^etre non nul. On sait alors avec \cite{image} 5.3.1 que l'image de cet op\'erateur est une repr\'esentation irr\'eductible. Comme par d\'efinition:
$$N^L(\rho_{v},\pi_{v},s)=\phi(\rho_{v},\pi_{v},s)N^{\GL}(\rho_{v},\pi_{v},s),$$ cet op\'erateur est lui aussi   holomorphe en $s=s_{0}$ et son image en $s=s_{0}$ contient, d'apr\`es ce que l'on a montr\'e, le sous-quotient de Langlands de l'induite $\rho_{v}\vert\,\vert^{s_{0}}\times \pi_{v}$. Par irr\'eductibilit\'e, ce sous-quotient est l'image de $N^{\GL}(\rho_{v},\pi_{v},s_{0})$ comme annonc\'e dans la proposition.

R\'eciproquement, supposons que $N^{\GL}(\rho_{v},\pi_{v},s_{0})$ est d'image non nulle \'egale au sous-quotient de Langlands de l'induite $\rho_{v}\vert\,\vert^{s_{0}}\times \pi_{v}$, que l'on note ici $\pi_{0}$. Or l'op\'erateur $N(\sigma,s_{0})$ est non nul donc son image contient $\pi_{0}$ qui est l'unique sous-module irr\'eductible de $\tau(s_{0})^*\times \pi_{temp}$. Mais $\pi_{0}$ intervient avec multiplicit\'e au plus un dans toutes les induites \'ecrites et ainsi $N(\sigma,s_{0})$ ne peut \^etre nul sur $\pi_{0}$ qui est l'image de $N^{\GL}(\rho_{v},\pi_{v},s_{0})$ par hypoth\`ese. De l'\'egalit\'e
$0\neq$
$$
N(\sigma,s_{0})\circ N^{\GL}(\rho_{v},\pi_{v},s_{0})=\phi(\rho_{v},\pi_{v},s_{0}) \bigl(N(\sigma,s)\circ N^{\GL}(\rho_{v},\pi_{v},s)\bigl)_{s=s_{0}},
$$
on obtient que $\phi(\rho_{v},\pi_{v},s_{0})\neq 0$. Ce qui termine la preuve.

\begin{cor} Pour tout $s_{0}\in {\mathbb C}$, la fonction $L(\rho_{v}\times \pi_{v},s)/L(\rho_{v}\times \pi^{\GL}_{v},s)$ est holomorphe en $s=s_{0}$.
\end{cor}
On montre d'abord le corrolaire en tout point $s_{0}\in {\mathbb R}_{\geq 1/2}$. Fixons $a\in {\mathbb N}$ un entier tr\`es grand. Sur les d\'efinitions, on v\'erifie que ni $L_{v}(\rho_{v}\times \pi_{v},s+a)$ ni $L_{v}(\rho_{v}\times \pi^{{\GL}}_{v},s+a)$ n'ont de p\^oles en $s=s_{0}$ (et il n'y a jamais de z\'eros). On reprend la notation $\phi(\rho_{v},\pi_{v},s)$ d\'ej\`a introduite. On a \'evidemment pour tout $a\in {\mathbb N}$:
$$
L(\rho_{v},\pi_{v},s)/L(\rho_{v}\times \pi^{\GL}_{v},s)=L(\rho_{v},\pi_{v},s+a)/L(\rho_{v},\pi^{\GL}_{v},s+a)$$
$$\prod_{i\in [0,a[}\phi(\rho_{v},\pi_{v},s+i).
$$
Donc l'holomorphie cherch\'ee est bien un corollaire de la proposition pr\'ec\'edente. On d\'emontre maitenant l'assertion pour $s_{0}\in {\mathbb R}_{<1/2}$. On utilise l'\'egalit\'e des facteurs $\gamma$ d\'emontr\'e en \ref{gamma}:
$$
L(\rho_{v}\times \pi_{v},s)/L(\rho_{v}\times \pi^{GL}_{v},s)=\zeta(s)L(\rho_{v}^*\times \pi_{v},1-s)/L(\rho_{v}^*\times \pi^{\GL}_{v},1-s),
$$
o\`u $\zeta(s)$ prend en compte les facteurs $\epsilon$ et est donc une fonction holomorphe inversible. Le terme de droite est holomorphe en $s_{0}$ puisque $1-s_{0}\in {\mathbb R}_{>1/2}$ et le terme de gauche est donc aussi holomorphe en $s=s_{0}$. Ensuite, on passe de $s_{0}$ r\'eel \`a $s_{0}$ complexe en tensorisant $\rho_{v}$ par un caract\`ere unitaire.
\subsection{P\^oles des facteurs $L$, le cas des places archim\'e\-diennes\label{localarchimedien}}
On s'attend bien \'evidemment \`a avoir les m\^emes r\'esultats aux places archim\'ediennes que ceux d\'emntr\'es au paragraphe pr\'ec\'edent pour les places finies. Toutefois, dans l'\'etat actuel des connaissances, il nous faut mettre des hypoth\`eses pour obtenir ces r\'esultats. Soit ici $v$ une place archim\'edienne de $k$.

On fixe donc $\pi$ et $\rho$ comme pr\'ec\'edemment et on fixe un nombre r\'eel positif $s_{1}$; d'o\`u la repr\'esentation $\pi^{\GL}$ de $\GL(m^*_{G},k_{v})$. L'hypoth\`ese forte faite ici est que $\pi_{0}$ a de la cohomologie et le caract\`ere infinit\'esimal de l'induite $\rho_{v}\vert\,\vert^{s_{1}}\times \pi_{v}$ est entier et r\'egulier; entier veut dire que si l'on identifie le caract\`ere infinit\'esimal \`a une collection de demi-entiers, ces demi-entiers sont soient tous entiers soit tous demi-entier non entier.

On fixe $s_{0}$ un nombre r\'eel et on suppose que $s_{0}\geq s_{1}$.
\begin{prop}(i) $\ord_{s=s_{0}}\bigl(L_{v}(\rho_{v}\times \pi_{v},s)/L_{v}(\rho_{v}\times \pi_{v},s+1)\bigr)\geq$
$$\ord_{s=s_{0}}\bigl(L_{v}(\rho_{v}\times \pi^{\GL}_{v},s)/L_{v}(\rho_{v}\times \pi^{\GL}_{v},s+1)\bigr).$$
(ii) $\ord_{s=s_{0}}L_{v}(\rho_{v}\times \pi_{v},s) \geq \ord_{s=s_{0}}L_{v}(\rho_{v}\times \pi^{\GL}_{v},s).$
\end{prop}
On d\'emontre un r\'esultat beaucoup plus fort qui lui n'a aucune chance d'\^etre vrai sans l'hypoth\`ese sur $s_{1}$: 
$$\ord_{s=s_{0}}L_{v}(\rho_{v}\times \pi_{v},s) =\ord_{s=s_{0}}L_{v}(\rho_{v}\times \pi^{\GL}_{v},s)=0.$$
Une fois ce r\'esultat d\'emontr\'e, on d\'eduit que toute les fonctions intervenant dans la proposition sont holomorphes non nulles en $s=s_{0}$, ce qui montre a fortiori la proposition.

Avant de faire la preuve remarquons que l'hypoth\`ese sur $s_{0}$ est un peu plus g\'en\'erale que de supposer que l'induite $\rho_{v}\vert\,\vert^{s_{0}}\times \pi_{v}$ a un caract\`ere infinit\'esimal r\'egulier.

On exploite le fait que $\pi_{v}$ a de la cohomologie: ainsi on sait que $\pi_{v}$ est dans l'un des paquet d'Adams-Johnson. On  note $\pi^{AJ}_{v}$ la repr\'esentation de $\GL(m^*_{G},k_{v})$ qui param\'etrise le paquet d'Adams-Johnson que l'on fixe arbitrairement. On \'ecrit $\pi^{AJ}_{v}$ sous la forme d'une induite de module de Speh,$$
\times_{\delta} {{\Speh}}(\delta,b_{\delta}) \times {{\Speh}}(\eta,b_{\eta})
$$o\`u $\delta$ parcourt un ensemble (sans multiplicit\'e) de s\'eries discr\`etes et o\`u $\eta$ est un caract\`ere quadratique,
auxquels s'ajoute encore \'eventuellement un caract\`ere quadratique $\eta'$ de ${\mathbb R}^*$ (cf\ref{gamma}). Comme dans \ref{gamma}, on \'ecrit $\pi_{v}$ comme (sous)-quotient de Langlands d'une induite de la forme
$$
\times_{\delta}\times_{ i\in [(b_{\delta}-1)/2,b_{\delta}'+1]}\delta\vert\,\vert^{i}\times_{j\in [(b_{(\eta-1)/2,(b'_{\delta}+1)/2]})}\times \pi_{temp},
$$
o\`u $\pi_{temp}$ est une repr\'esentation temp\'er\'ee convenable; il faudrait r\'eor\-donner les repr\'esentations pour que $\pi_{v}$ soit r\'eellement un quotient.
Il suffit de montrer que $L_{v}(\rho_{v}\times \delta,s-i)$ est holomorphe en $s=s_{0}$ pour tout $\delta$, $i$ comme ci-dessus et qu'il en est de m\^eme pour $L_{v}(\rho_{v}\times \eta,s-j)$ avec les notations pr\'ec\'edentes. On peut supposer que $\rho_{v}$ est soit une s\'erie discr\`ete soit un caract\`ere unitaire: en g\'en\'eral on sait que $\rho_{v}$ est le produit de telles repr\'esentations par des caract\`eres de la forme $\vert\,\vert^x$ o\`u $x\in ]-1/2,1/2[$. Pour un tel $x\neq 0$ l'holomorphie cherch\'ee r\'esulte simplement de ce que $\pi_{v}$ a un caract\`ere infinit\'esimal entier. Le cas des caract\`eres unitaires non quadratiques se r\`egle de fa\c{c}on analogue. On suppose donc que $x=0$ on commence par le cas des s\'eries discr\`etes: soit $\delta'$ une s\'erie discr\`ete de $\GL(2,{\mathbb R})$ et on note $r'$ son param\`etre et $t'=r'/2$, c'est-\`a-dire que le caract\`ere infinit\'esimal de $\delta'$ est $t',-t'$. On calcule les facteurs $L$ de paires en utilisant les formules \cite{shahidiarchimedien} (3.5) et (3.7) qui n\'ecessitent de savoir d\'ecomposer le produit tensoriel de deux repr\'esentations du groupe de Weil de ${\mathbb R}$; cette d\'ecomposition est particuli\`erement bien \'ecrite dans \cite{kasten} proposition 1.6. On a donc \`a consid\'erer des p\^oles de fonctions de l'une des formes suivantes:

pour $\delta$ intervenant ci-dessus, de param\`etre not\'e $t$ et pour $\ell\in [(b_{\delta}-1)/2,-(b_{\delta}-1)/2]$,
 $\Gamma(s-\ell+\vert t'-t\vert)$;

 et pour  $\ell' \in [(b_{\eta}-1)/2,-(b_{\eta}-1)/2]$, $\Gamma (s+t-\ell')$.

Consid\'erons le premier type de fonctions; elles ont des p\^oles uniquement si $s-(b_{\delta}-1)/2+\vert t-t'\vert\leq 0$ et si ce nombre est entier relatif. Ce qui se d\'eveloppe en :
$$
t'+s\leq (b_{\delta}-1)/2+t \hbox{ et } -t'+s\leq -t+(b_{\delta}-1)/2.
$$
Mais le caract\`ere infinit\'esimal de la repr\'esentation ${\Speh}(\delta,b_{\delta})$ est la collection de demi-entiers $ t+[(b_{\delta}-1)/2,-(b_{\delta}-1)/2] \cup -t+ [(b_{\delta}-1)/2,-(b_{\delta}-1)/2]$. Par hypoth\`ese le caract\`ere infinit\'esimal de $\delta'\vert\,\vert^{s_{1}}\times {\Speh}(\delta,b_{\delta})$ est certainement entier et r\'egulier. On a donc soit $t'+s_{1}>t+(b_{\delta}-1)/2$ soit 
$$
t-(b_{\delta}-1)/2>t'+s_{1}\geq -t'+s_{1}>-t+(b_{\delta}-1)/2.
$$
Puisque $s_{0}\geq s_{1}$, on a soit $t'+s_{0}>t+(b_{\delta}-1)/2$ soit $-t'+s_{0}>-t+(b_{\delta}-1)/2$ et dans tous les cas, les conditions pour avoir un p\^ole ne sont pas satisfaites. 

Comme le caract\`ere infinit\'esimal de $\delta\vert\,\vert^{s_{1}}\times {\Speh}(\eta,b_{\eta})$ et lui aussi entier et r\'egulier, on v\'erifie encore plus facilement que le deuxi\`eme type de fonctions n'ont pas de p\^ole en $s=s_{0}$. Si $\delta'$ est remplac\'e par un caract\`ere, on raisonne de la m\^eme fa\c{c}on.

Ici on a en fait d\'emontrer, au passage, que la fonction $L(\rho_{v}\times \pi^{\GL}_{v},s)$ n'a pas de p\^ole en $s=s_{0}$ puisque $\pi^{\GL}_{v}$ a une forme de m\^eme type que $\pi^{AJ}$ et a le m\^eme caract\`ere infinit\'esimal; pour les fonctions $L(\rho_{v}\times \pi_{v},s)$ il faut encore  remarquer que par positivit\'e stricte la fonction $L(\rho_{v}\times \pi_{temp},s)$ n'a pas de p\^ole en $s=s_{0}$ et cela termine la d\'emonstration.
\subsection{Op\'erateurs d'entrelacement aux places archim\'e\-diennes\label{entrelacementarchimedien}}
On continue avec les hypoth\`eses du paragraphe pr\'ec\'edents: $v$ est une place archim\'edienne et $\pi_{v}$ a de la cohomologie. On suppose qu'il existe un nombre r\'eel $s_{1}>0$ tel que le caract\`ere infinit\'esimal de l'induite $\rho_{v}\vert\,\vert^{s_{1}}\times \pi_{v}$ soit entier et r\'egulier.
\begin{lem} Pour tout nombre r\'eel $s_{0}\geq s_{1}$ l'induite $\rho_{v}\vert\,\vert^{s_{1}}\times \pi_{v}$ a un unique quotient irr\'eductible qui est le sous-quotient de Langlands. L'o\-p\'e\-rateur d'entrelacement standard qui envoie l'induite $\rho_{v}\vert\,\vert^{s}\times \pi_{v}$ dans $\rho^*_{v}\vert\,\vert^{-s}\times \pi_{v}$ est holomorphe en $s=s_{0}$ et son image est une repr\'esenta\-tion irr\'eductible isomorphe au sous-quotient de Langlands de l'induite.
\end{lem}
On \'ecrit $\pi_{v}$ comme quotient de Langlands d'une induite $\times_{i\in [1,\ell]}\tau_{i}\times \pi_{temp}$, o\`u les $\tau_{i}$ sont des s\'eries discr\`etes ou des caract\`eres unitaires tordus dont les exposants sont positifs strictement et rang\'es dans l'ordre d\'ecroissant. On peut aussi supposer que $\pi_{v}$ est une s\'erie discr\`ete ou un caract\`ere quadratique car en $s=s_{1}$, on a suppos\'e que l'induite $\rho_{v}\vert\,\vert^{s_{1}}\times \pi_{v}$ a un caract\`ere entier. On fixe $s_{0}$ comme dans l'\'enonc\'e et  le point est que si $s_{0}$ est plus petit (au sens large) que l'exposant de $\tau_{i}$ pour $i\in [1,\ell]$ fix\'e, alors l'induite $\rho_{v}\vert\,\vert^{s_{0}}\times \tau_{i}$ d'un groupe lin\'eaire convenable est irr\'eductible; cela r\'esulte de la comparaison des param\`etres de $\rho_{v}\vert\,\vert^{s_{0}}$ et de $\tau_{i}$ comme on l'a fait dans le paragraphe pr\'ec\'edent. Ainsi l'induite $\rho_{v}\vert\,\vert^{s_{0}}\times \pi_{v}$ est un quotient de l'induite $\rho_{v}\vert\,\vert^{s_{0}}\times_{i\in [1,\ell]}\tau_{i}\times \pi_{temp}$ qui est elle-m\^eme isomorphe \`a l'induite r\'eordonn\'ee pour que les exposants soient rang\'es dans un ordre d\'ecroissant. Cela prouve la premi\`ere partie de l'\'enonc\'e.

Il reste \`a prouver l'holomorphie de l'op\'erateur d'entrelacement.
On v\'erifie aussi que sous ces m\^emes hypoth\`eses, l'op\'erateur d'entrelacement standard est holomorphe: on le sait si $s_{0}$ est strictement plus grand que l'exposant de $\tau_{i}$. Si les exposants sont \'egaux, on sait que l'op\'erateur d'entrelacement normalis\'e \`a la Langlands-Shahidi est holomorphe mais comme on a montr\'e que le facteur de normalisation n'a ni z\'ero ni p\^ole, on obtient l'assertion.
\begin{rmq} Les propri\'et\'es de l'op\'erateur d'entrelacement standard du lemme pr\'ec\'edent se transf\`erent imm\'ediatement \`a l'op\'erateur d'entrelacement normalis\'e,  $N^{\GL}(\rho_{v}\times \pi_{v},s)$ d\'efini comme en \ref{localfini} puisque le facteur de normalisation n'a ni z\'ero ni p\^ole.
\end{rmq}
\subsection{Remarque sur les places archim\'ediennes\label{remarquearchimedien}}
On a d\'emontr\'e l'\'equation fonctionnelles aux places archim\'ediennes sous l'hypoth\`ese que $\pi^{\GL}_{v}$ est le param\`etre d'un paquet d'Adams-John\-son contenant $\pi_{v}$.
\begin{rmq} La proposition \ref{localfini} est aussi vraie \`a la place archim\'edienne $v$ et m\^eme plus pr\'ecis\'ement ici (i) est vrai pour tout $s_{0}\in {\mathbb R}_{>0}$.
\end{rmq}
Cette remarque est, en fait, compl\`etement triviale: on d\'ecrit les param\`etres de Langlands de $\pi_{v}$ en fonction de $\pi^{\GL}_{v}$ (au moins partiellement) comme on la fait en \ref{gammaarchimedien}. D'apr\`es la d\'efinition des facteurs $L$ qui est multiplicative en les param\`etres de Langlands, on est directement ramen\'e au cas o\`u $\pi_{v}$ est temp\'er\'e. On suppose d'abord que $s_{0}\in {\mathbb R}_{>0}$. Par positivit\'e, $L_{v}(\rho_{v}\times \pi_{v},s)$ n'a pas de p\^ole en $s=s_{0}$ et n'y est pas nul car il ne peut y avoir de $0$. On a donc de fa\c{c}on imm\'ediate le (ii) de \ref{localfini} sous l'hypoth\`ese que $s_{0}>0$; gr\^ace \`a l'\'equation fonctionnelle, on obtient comme dans \ref{localfini} l'assertion (ii) pour tout $s_{0}$. Pour montrer (i), on calcule avec la description explicite de $\pi^{\GL}_{v}$: $$L(\rho_{v}\times \pi^{GL}_{v},s)/L(\rho_{v}\times \pi^{GL}_{v},s+1)=\prod_{(\delta,b_{\delta})}\prod_{k\in [(b_{\delta}-1)/2,-(b_{\delta}-1)/2]}
$$
$$L(\rho_{v}\times \delta,s-k)/L(\rho_{v}\times \sigma,s-k+1),
$$
auquel s'ajoute un produit du m\^eme type o\`u $\delta$ est remplac\'e par un  caract\`ere. Mais un tel produit se simplifie pour donner
$$
\prod_{\delta,b_{\delta}}L(\rho_{v}\times \delta,s-(b_{\delta}-1)/2)/L(\rho_{v}\times \delta,s+(b_{\delta}+1)/2).
$$
Les d\'enominateurs sont donc d'ordre $0$ en $s=s_{0}\in {\mathbb R}_{>0}$ et en un tel point l'ordre du quotient de fonctions $L$ est donc inf\'erieur ou \'egal \`a $0$ ce qui donne (i).

\subsection{P\^oles des fonctions $L$\label{poleglobal}}
On fixe $\rho$ et $\pi$ et on reprend les hypoth\`eses des places archim\'edien\-nes: on suppose qu'il existe un r\'eel positif $s_{1}$ tel qu'en toute place archim\'edienne $v$ le caract\`ere infinit\'esimal de l'induite  $\rho_{v}\vert\,\vert^{s_{1}}\times \pi_{v}$ est entier et r\'egulier.

\begin{thm} (i) Soit $s_{0}\in {\mathbb R}_{\geq s_{1}}$.
$\ord_{s=s_{0}}L(\rho\times \pi,s)/L(\rho\times \pi,s+1)\geq$
$$
\ord_{s=s_{0}}L(\rho\times \pi^{\GL},s)/L(\rho\times \pi^{\GL},s+1).$$
et $\ord_{s=s_{0}}L(\rho\times \pi,s) \geq \ord_{s=s_{0}}L(\rho\times \pi^{\GL},s)$.

(ii) On suppose ici qu'en toute place archim\'edienne, r\'eelle, $v$, $\pi^{\GL}_{v}$ param\'etrise un paquet d'Adams-Johnson contenant $\pi_{v}$. Alors $L(\rho\times \pi,s)$ et $L(\rho\times \pi^{\GL},s)$ v\'erifient la m\^eme \'equation fonctionnelle et pour tout $s_{0}\in {\mathbb C}$:
$$
\ord_{s=s_{0}}L(\rho\times \pi,s)/L(\rho\times \pi,s+1)\geq
\ord_{s=s_{0}}L(\rho\times \pi^{\GL},s)/L(\rho\times \pi^{\GL},s+1).$$
$$
\ord_{s=s_{0}}L(\rho\times \pi,s)\geq \ord_{s=s_{0}}L(\rho\times \pi^{\GL},s).
$$

\end{thm}
Le (i) de ce th\'eor\`eme est un corollaire des r\'esultats locaux. Ceci est aussi vrai pour le (ii), pour l'\'equation fonctionnelle et pour la comparaison des ordres si l'on suppose que $s_{0}\in {\mathbb R}_{\geq 1/2}$ puisque l'on a cette limitation pour les places finies.
Soit donc $s_{0}\in {\mathbb R}_{\leq 1/2}$; alors $1-s_{0}\in {\mathbb R}_{\geq 1/2}$ et par l'\'equation fonctionnelle, en posant $s'=1-s$
$$
L(\rho^*\times \pi,s)/L(\rho^*\times \pi,s+1)=\epsilon(\rho^*\times \pi,s)L(\rho\times \pi,s')/L(\rho\times \pi,1+s').
$$
L'ordre du membre de gauche en $s=s_{0}$ est donc \'egal \`a celui du membre de droite en $s'=1-s_{0}$. On a la m\^eme \'egalit\'e en rempla\c{c}ant $\pi$ par $\pi^{\GL}$ et cela permet de d\'emontrer (ii) pour tout $s_{0}$.

\section{Applications aux s\'eries d'Eisenstein}
Le but de cette section est de d\'emontrer que l'existence d'un p\^ole en $s=s_{0}\in {\mathbb R}_{> 1/2}$ pour la fonction $L(\rho\times \pi,s)/L(\rho\times \pi,s+1)$ entra\^{\i}ne l'existence d'un p\^ole pour les s\'eries d'Eisenstein $E(\rho\times \pi,s)$ en $s=s_{0}$. On donne m\^eme une condition n\'ecessaire et suffisante. Le cas o\`u $s=1/2$ et un peu diff\'erent, on sait par avance que les p\^oles des s\'eries d'Eisenstein en $s=1/2$ d\'ependent non pas de la fonction pr\'ec\'edente (dont on v\'erifiera qu'elle est holomorphe) mais d'une fonction $L$ li\'ee \`a la repr\'esentation $\rho$ elle-m\^eme, la fonction $L(\rho,r_{G},s)$ o\`u $r_{G}$ est la repr\'esentation naturelle du groupe dual de $G$ dans $GL(m^*_{G},{\mathbb C})$. Et en $s_{0}=1/2$ ce sont les p\^oles de la fonction $L(\rho,r_{G},s)L(\rho\times \pi,s)/L(\rho\times \pi,s+1)$ qui vont intervenir.

Malheureusement il nous faut des hypoth\`eses aux places archim\'e\-dien\-nes: on suppose qu'il existe un r\'eel $s_{1}$ strictement positif tel qu'en toute place archim\'edienne $v$,  l'induite $\rho_{v}\vert\,\vert^{s_{1}}\times \pi_{v}$ a un caract\`ere infinit\'esimal entier et r\'egulier.

\subsection{P\^oles des quotients de fonctions L et r\'esidus de s\'eries d'Eisenstein\label{eisenstein}}
On reprend la d\'efiniton du sous-quotient de Langlands des induites locales $\rho_{v}\vert\,\vert^{s_{0}}\times \pi_{v}$,  d\'ej\`a utilis\'e en \ref{localfini}, que l'on note ici $Lang(\rho_{v}\vert\,\vert^{s_{0}}\times \pi_{v})$.
\begin{thm} On fixe un nombre r\'eel $s_{0}\in {\mathbb R}_{\geq s_{1}}$. Et on suppose que $\pi$ est une repr\'esentation cuspidale.

(i) On suppose ici que $s_{0}>1/2$; alors le p\^ole de la fonction $L(\rho\times \pi,s)/L(\rho\times \pi,s+1)$ est au plus simple en $s=s_{0}$. Et on a l'\'equivalence:
$L(\rho\times \pi,s)/L(\rho\times \pi,s+1)$ a un p\^ole en $s=s_{0}$  si et seulement si la s\'erie d'Eisenstein $E(\rho\times \pi,s)$ a un p\^ole en $s=s_{0}$ et son r\'esidu est une repr\'esentation n\'ecessairement irr\'eductible isomorphe en toute place $v$ \`a $Lang(\rho_{v}\vert\,\vert^{s_{0}}\times \pi_{v})$.

(ii)On suppose ici que $s_{0}=1/2$; alors la fonction $L(\rho\times \pi,s)/L(\rho\times \pi,s+1)$ est holomorphe en $s=s_{0}$ et l'on a l'\'equivalence

$
L(\rho,r_{G},2s)L(\rho\times \pi,s)/L(\rho\times \pi,s+1)$ a un p\^ole en $s=s_{0}$ si et seulement si $$
\bigl((s-s_{0})E(\rho\times \pi,s)\bigr)_{s=s_{0}}\simeq \otimes'_{v}Lang(\rho_{v}\vert\,\vert^{s_{0}}\times \pi_{v}).
$$
\end{thm}
(i) et (ii) se d\'emontrent simultan\'ement. Puisque l'on a suppos\'e que $\pi$ est cuspidale,
les s\'eries d'Eisenstein consid\'e\-r\'ees n'ont qu'un terme constant non nul qui vaut, pour $f_{s}$ une section de l'induite $\rho\vert\,\vert^{s}\times \pi$,
$
f_{s}+M(s)f_{s}
$,
o\`u $M(s)$ est l'op\'erateur d'entrelacement standard. On remplace $M(s)$ par $$\frac{L(\rho, r_{G},2s)}{L(\rho,r_{G},2s+1)}\frac{L(\rho\times \pi^{{\GL}},s)}{L(\rho\times \pi^{{\GL}},s+1)} N^{{\GL}}(\rho\times \pi,s),$$ o\`u $N^{{\GL}}(s)$ est le produit sur toutes les places $v$ de l'op\'erateur local $N_{v}^{{\GL}}(\rho_{v}\times \pi_{v},s)$ d\'ej\`a utilis\'e en \ref{localfini} a un facteur pr\`es qui vient des facteurs locaux de $L(\rho, r_{G},2s)/L(\rho,r_{G},2s+1)$. Ces facteurs locaux n'ont ni z\'ero  ni p\^ole en $s=s_{0}\geq  1/2$. On a donc
$
\ord_{s=s_{0}}M(s)=\ord_{s=s_{0}}L(\rho, r_{G},2s)/L(\rho,r_{G},2s+1)$
$$+\, \ord_{s=s_{0}}L(\rho\times \pi^{{\GL}},s)/L(\rho\times \pi^{{\GL}},s+1) $$
$$+\sum_{v}\ord_{s=s_{0}} N_{v}^{{\GL}}(\rho_{v}\times \pi_{v},s),
$$
la somme sur $v$ \'etant n\'ecessairement finie. Et ceci vaut encore, avec les notations de \ref{localfini}
$$\ord_{s=s_{0}}L(\rho, r_{G},2s)/L(\rho,r_{G},2s+1)+
ordre_{s=s_{0}}L(\rho\times \pi,s)/L(\rho\times \pi,s+1)$$
$$+\sum_{v} \ord_{s=s_{0}}\phi(\rho_{v},\pi_{v},s)^{-1}+ordre_{s=s_{0}} N_{v}^{{\GL}}(\rho_{v}\times \pi_{v},s).
$$Dans la somme sur les places $v$, on est en droit de se limiter aux places finies d'apr\`es \ref{localarchimedien}.
Le lemme \ref{localfini} montre que chaque terme de cette somme est inf\'erieur ou \'egal \`a $0$. Et on a donc
$$
\ord_{s=s_{0}}M(s)\leq \ord_{s=s_{0}}L(\rho, r_{G},2s)/L(\rho,r_{G},2s+1)$$
$$+
 \ord_{s=s_{0}}L(\rho\times \pi,s)/L(\rho\times \pi,s+1).
$$
On rappelle que Langlands a d\'emontr\'e que ces s\'eries d'Eisenstein ont des p\^oles au plus simple  \cite{langlands} et que les p\^oles de $M(s)$ sont exactement les p\^oles de ces s\'eries d'Eisenstein. Donc on retrouve d\'ej\`a que le produit
$$
\bigl(L(\rho, r_{G},2s)/L(\rho,r_{G},2s+1)\bigr)\bigl(L(\rho\times \pi^{{\GL}},s)/L(\rho\times \pi^{{\GL}},s+1)\bigr) 
$$
a au plus des p\^oles simples.

De plus, puisque $s_{0}\geq 1/2$ les p\^oles de $L(\rho, r_{G},2s)/L(\rho,r_{G},2s+1)$ sont parmi ceux de $L(\rho\times \rho,2s)$; il y en a donc au plus un, n\'ecessairement simple, en $s_{0}=1/2$ et  $L(\rho, r_{G},2s)/L(\rho,r_{G},2s+1)$ n'a certainement pas de $0$ pour $s_{0}>1/2$ puisque cela fournirait un $0$ hors de la bande critique (ouverte) \`a $L(\rho\times \rho,2s)$. De plus, puisque les p\^oles de $L(\rho\times \pi,s)/L(\rho\times \pi,s+1)$ fournissent des p\^oles \`a $L(\rho\times \pi^{GL},s)/L(\rho\times \pi^{GL},s+1)$ et qu'il n'y en a donc pas en $s=1/2$. Cela montre les propri\'et\'es d'holomophie ou d'ordre des p\^oles annonc\'ees dans l'\'enonc\'e concernant cette fonction.

On en vient au c{\oe}ur de la d\'emonstration: on a une \'equivalence entre le fait que $$\bigl(L(\rho, r_{G},2s)/L(\rho,r_{G},2s+1)\bigr)\bigl(L(\rho\times \pi^{{\GL}},s)/L(\rho\times \pi^{{\GL}},s+1)\bigr) $$ ait un p\^ole en $s=s_{0}$ et le fait que d'une part les s\'eries d'Eisenstein $E(\rho\times\pi,s)$ aient un p\^ole en $s=s_{0}$ et d'autre part que pour toute place $v$:
$$
 \ord_{s=s_{0}}\phi(\rho_{v},\pi_{v},s)=ordre_{s=s_{0}} N_{v}^{{\GL}}(\rho_{v}\times \pi_{v},s)\eqno(1)$$
 Avec les hypoth\`eses que l'on a mises, si $v$ est une place archim\'edienne, on sait que (1) est vrai et que l'image de $N^{\GL}(\rho_{v}\times \pi_{v},s_{0})$ est le sous-quotient de Langlands. Aux places finies, on remplace (1) par la propri\'et\'e \'equivalente, d'apr\`es \ref{localfini}, que le sous-quotient de Langlands de l'induite $\rho_{v}\vert\,\vert^{s_{0}}\times \pi_{v}$ est exactement l'image de $N_{v}^{\GL}(\rho_{v}\times \pi_{v},s_{0})$. Cela donne exactement les \'equivalences en (i) et (ii) de l'\'enonc\'e du th\'eor\`eme.
\subsection{Images par s\'eries theta \label{theta}} 
Dans cette sous-section, on suppose que $G$ est un groupe sp\'ecial orthogonal, $SO(m)$ et on fixe $r$ un entier satisfaisant $2r<m$. On note $H$ le groupe $Sp(2r)$ si $m$ est pair et $Mp(2r)$ si $m$ est impair. On dit que $\pi$ est dans l'image cuspidale des s\'eries theta \`a partir de $H$ s'il existe une repr\'esentation cuspidale $\tau$ de $H$ et un caract\`ere quadratique automorphe $\eta$ telle que $\pi\otimes \eta$ soit inclus dans la restriction \`a $SO(m)$ de l'image de $\tau$ par s\'erie thetas. Comme il y a plusieurs fa\c{c}ons de normaliser l'image par s\'eries theta, puisque l'on peut toujours tordre par un caract\`ere automorphe de $SO(m)$, la d\'efinition que l'on a prise \'evite de devoir pr\'eciser. Comme cela, on a aussi une ind\'ependance en fonction du caract\`ere additif fix\'e pour d\'efinir les s\'eries theta.

On suppose que pour toute place archim\'edienne $v$, $\pi_{v}$ a de la cohomologie et qu'il existe un caract\`ere quadratique automorphe $\eta$ tel que l'induite $\eta_{v}\vert\,\vert^{m/2-r}\times \pi_{v}$ est entier r\'egulier. Cette derni\`ere propri\'et\'e est ind\'ependante du caract\`ere quadratique.
\begin{thm} On suppose que la fonction $L(\eta\times \pi,s)/L(\eta\times \pi,s+1)$ a un p\^ole en $s=m/2-r$ pour un choix de caract\`ere quadratique $\eta$; alors $\pi$ est dans l'image cuspidale des s\'eries theta \`a partir de $H$.
\end{thm}
D'apr\`es le th\'eor\`eme de \ref{eisenstein} dont les hypoth\`eses sont satisfaites, on a en fait n\'ecessairement $m/2-r\geq 1$ et les s\'eries d'Eisenstein $E(\eta\times \pi,s)$ ont un p\^ole en $s=m/2-r$. Il r\'esulte de \cite{lie} thm p203 et \cite{jiangsoudry} 5.1 qu'il existe un entier $r'\leq r$ tel que $\pi$ soit dans l'image cuspidale des s\'eries theta de $H'$ o\`u $H'$ est l'analogue de $H$ en rempla\c{c}ant $r$ par $r'$. On note $\tau'$ une repr\'esentation cuspidale et $\eta'$ un caract\`ere quadratique tel que  $\pi\otimes \eta'$ soit  dans l'image du rel\`evement theta de $\tau'$. En regardans la correspondance theta non ramifi\'e, on voit qu'il existe un caract\`ere $\eta''$ quadratique (qui d\'epend des choix) tel que $\pi^{\GL}$ soit une induite de module de Speh l'un d'entre eux \'etant n\'ecessairement ${\Speh}(\eta'', m-1-2r')$ et, uniquement si $m$ est pair, un facteur ${\Speh}(\eta''',1)$ o\`u $\eta'''$ est un caract\`ere quadratique convenable. Puisque $L(\eta\times \pi,s)/L(\eta\times \pi,s+1)$ a un p\^ole en $s=m/2-r$ il en est de m\^eme de $L(\eta\times \pi^{\GL},s)/L(\eta\times \pi^{\GL},s+1)$ d'apr\`es \ref{poleglobal}. Cela force $\pi^{\GL}$ \'ecrit comme induite de module de Speh d'avoir aussi un facteur ${\Speh}(\eta,m-1-2r)$. Si $r\neq r'$ la r\'egularit\'e du caract\`ere infinit\'esimal de $\pi$ aux places archim\'ediennes n'est s\^urement pas satisfaite. Ainsi $r'=r$ et le th\'eor\`eme est d\'emontr\'e.

\

A l'instar de \ref{eisenstein} on pourrait donner un \'enonc\'e o\`u l'on a une \'equivalence et non seulement une condition suffisante.

\subsection{Exemple\label{exemple}}
On propose ici un exemple d'un triplet $(\rho,\pi,s_{0})$ o\`u les s\'eries d'Eisensteine $E(\rho\times \pi,s)$ ont un p\^ole en $s=s_{0}$ alors que la fonction $L(\rho\times \pi,s)/L(\rho\times \pi,s+1)$ est holomorphe en $s=s_{0}$. Dans tout ce qui suit $tr_{?}$ signifie la repr\'esentation triviale du groupe en indice.

On fixe $\rho'$ une repr\'esentation cuspidale autoduale de ${\GL}(6)$; on suppose qu'il existe deux places finies, not\'ee $v_{i}$, pour $i=1,2$ telles que  $\rho_{v_{i}}$ est l'induite   du caract\`ere trivial et de $St(\eta,5)$, la repr\'esentation de Steinberg de ${\GL}(5,k_{v_{i}})$ tordue par le caract\`ere quadratique non ramifi\'e non trivial. Il n'y a aucune raison pour qu'une telle repr\'esentation n'existe pas mais je n'ai pas non plus l'assurance qu'elle existe car les conditions locales donnent des repr\'esentations temp\'er\'ees et non des s\'eries discr\`etes. On suppose en plus qu'\`a l'infini la repr\'esentation a de la cohomologie et l\`a on a du choix.

On consid\`ere les donn\'ees globales, $(tr_{{\GL}(1)},1), (\rho',5)$ o\`u $tr_{{\GL}(1)} $ est la repr\'esentation trivial de ${\GL}(1)$. Ceci fournit un paquet d'Arthur de repr\'esentation de carr\'e int\'egrables pour $Sp(30,k)$, le groupe symplectique de rang 15.

Le caract\`ere d'Arthur est trivial; il n'y a que des repr\'esentations de type orthogonale.

Pour $i=1,2$, il existe une (en fait deux) repr\'esentation $\pi_{i}$  dans le paquet local d'Arthur qui v\'efie une inclusion
$$
\pi_{i}\hookrightarrow \vert\, \vert^{-2}\times \pi_{i,cusp},\eqno(1)
$$
o\`u $\pi_{i,cusp}$ est une repr\'esentation cuspidale de $Sp(28,k_{v_{i}})$: en effet, localament, le param\`etre d'Arthur est la somme de 3 repr\'esentation irr\'eductibles de $W_{k_{v_{i}}}\times SL(2,{\mathbb C})\times SL(2,{\mathbb C})$:
$$tr_{W_{k_{v_{i}}}\times SL(2,{\mathbb C})\times SL(2,{\mathbb C})}\oplus tr_{W_{k_{v_{i}}}\times SL(2,{\mathbb C})}\otimes sp_{5}\oplus \eta_{W_{k_{v_{i}}}}\otimes sp_{5}\otimes sp_{5},
$$o\`u $sp_{5}$ est la repr\'esentation irr\'eductible de dimension $5$ de $SL(2,{\mathbb C})$. On a donn\'e en \cite{contemporary} une description compl\`ete des repr\'esentations dans ce paquet (cf. la version simplifi\'ee \cite{holomorphie} 2.3).

\begin{rmq}
Soit $\pi$ une repr\'esentation dans le paquet d'Arthur global d\'efini par les param\`etres pr\'ec\'edent. On suppose que $\pi_{v_{1}}\simeq \pi_{1}$. Alors $\pi$ est cuspidale.
\end{rmq}
En effet si $\pi$ n'est pas cuspidale, elle v\'erifie \cite{manuscripta} 1.3 et en toute place $v$, $\pi_{v}$ est un quotient de $\rho'_{v}\vert\,\vert^{2}\times \pi'[v]$, o\`u $\pi'_{v}$ est une repr\'esentation convenable. Mais ceci n'est pas vrai en les places $v_{i}$ (pour $i=1,2$) \`a cause de (1), puisque $\rho'$ n'est pas un caract\`ere.

\begin{lem}Soit $\pi$ une repr\'esentation dans le paquet d'Arthur global d\'efini pr\'ec\'edemment; on suppose que $\pi_{v_{1}}\simeq \pi_{1}$. Le facteur $L_{v _{1}}(tr_{{\GL}(1)}\times \pi,s)/L_{v _{1}}(tr_{{\GL}(1)}\times \pi,s+1)$ a un z\'ero en $s=1$. Pour $i=1,2$ les facteurs $L_{v_{i} }(tr_{{\GL}(1)}\times \pi_{v_{i}}^{\GL},s)/L{v_{i} }(tr_{{\GL}(1)}\times \pi_{v_{i}}^{\GL},s+1)$ sont holomorphes non nuls en $s=1$.
\end{lem}
La premi\`ere assertion vient de (1) qui donne la param\'etrisation de Langlands de $\pi_{1}$. Elle est \'evidemment aussi vraie pour $i=2$. Montrons la deuxi\`eme; 
$$
\frac{L_{v_{i} }(tr_{{\GL}(1)}\times \pi_{v_{i}}^{\GL},s)}{L{v_{i} }(tr_{{\GL}(1)}\times \pi_{v_{i}}^{\GL},s+1)}=
\frac{L_{v_{i} }(s)}{L_{v_{i}}(s+1)}\frac{L_{v_{i}}(\rho'_{v_{i}}, s-2)}{L_{v_{i}}(\rho'_{v_{i}},s+3)}.
$$
La non nullit\'e est claire puisque les d\'enominateurs sont calcul\'es en des r\'eels strictement positifs. On doit donc montrer que pour $i=1,2$, $
L_{v_{i}}(s)L_{v_{i}}(\rho'_{v_{i}},s-2)$ n'a pas de p\^ole en $s=1$. Ceci est clair pour le premier facteur, le deuxi\`eme facteur est encore un produit, $L_{v_{i}}(s-2)L_{v_{i}}(St(\eta,5),s-2)$. Comme $\eta$ est non trivial, aucun des facteurs n'a de p\^ole en $s=1$.

\begin{rmq} Il existe une repr\'esentation cuspidale dans le paquet d'Arthur assoic\'ee aux donn\'ees pr\'ec\'edentes, telle que $\pi_{v_{1}}\simeq \pi_{1}$ et pour toute place $v\neq v_{1},v_{2}$, $\pi_{v}$ soit dans le paquet de Langlands associ\'ee au paquet d'Arthur. Pour une telle repr\'esentation $\pi$, les s\'eries d'Eisenstein $E(tr_{{\GL}(1)}\times \pi,s)$ ont un p\^ole en $s=1$ et le quotient de fonction $L$, $L(tr_{{\GL}(1)}\times \pi,s)/L(tr_{{\GL}(1)}\times \pi,s+1)$ est holomorphe en $s=1$.
\end{rmq}
Ce qui se produit est que des facteurs locaux pour le quotient de fonction $L$ ont des z\'eros alors qu'il n'y en a pas pour l'analogue avec ${\GL}$ en exposant mais l'op\'erateur d'entrelacement local n'est pas holomorphe quand il est normalis\'e \`a la Langlands alors qu'il l'est quand il est normalis\'e avec les facteurs venant du groupe g\'en\'eral lin\'eaire.

L'existence de $\pi$ vient de la formule de multiplicit\'e d'Arthur et on a gard\'e une libert\'e en la place $v_{2}$ pour assurer tout probl\`eme.

On a alors:
$$\frac{
L(tr_{{\GL}(1)}\times \pi,s)}{L(tr_{{\GL}(1)}\times \pi,s+1)}=\gamma_{v_{1}}\gamma_{v_{2}}
\frac{L(tr_{{\GL}(1)}\times \pi^{\GL},s)}{L(tr_{{\GL}(1)}\times \pi^{\GL},s+1)},
$$
o\`u pour $i=1,2$, 
$$
\gamma_{v_{i}}=\frac{L_{v_{i}}(tr_{{\GL}(1)}\times \pi_{v_{i}},s)}{L_{v_{i}}(tr_{{\GL}(1)}\times \pi_{v_{i}},s+1)}\frac {L_{v_{i}}(tr_{{\GL}(1)}\times \pi^{\GL}_{v_{i}},s+1)}{L_{v_{i}}(tr_{{\GL}(1)}\times \pi^{\GL}_{v_{i}},s)}.
$$
Ces facteurs n'ont pas de p\^oles pour $s\in {\mathbb R}_{>0}$ par l'argument donn\'e ci-dessus. Pour $i=1,2$, en $s=1$, $\gamma_{v_{i}}$ a un z\'ero d'apr\`es les r\'esultats pr\'ec\'edents. On sait que $L(tr_{{\GL}(1)}\times \pi^{{\GL}},s)/L(tr_{{\GL}(1)}\times \pi^{{\GL}},s+1)$ a un p\^ole d'ordre exactement 1 en $s=1$ \`a cause du premier facteur. Ainsi ce qui pr\'ec\`ede montre que $L(tr_{{\GL}(1)}\times \pi,s)/L(tr_{{\GL}(1)}\times \pi,s+1)$ est holomorphe en $s=1$.

Il reste \`a montrer que les s\'eries d'Eisenstein ont des r\'esidus non nul en $s=1$. Avec \cite{holomorphie}, il faut encore v\'erifier que les op\'erateurs d'entrelacements locaux normalis\'es comme en loc.cite sont non nuls. Mais cela r\'esulte des choix: aux places $v_{i}$, pour $i=1,2$, la non nullit\'e r\'esulte de ce que l'op\'erateur d'entrelacement standard:
$$
\vert\,\vert^{s}\times \vert\,\vert^{-2}\times \pi_{v_{i},cusp}\rightarrow \vert\,\vert^{-s}\times \vert\,\vert^{-2}\times \pi_{v_{i},cusp},
$$
est holomorphe et donc n\'ecessairement non nul; la normalisation avec les fonctions $L$ attach\'ees \`a $\pi^{\GL}$, n'introduit ni p\^ole ni z\'ero \`a cette place comme on l'a vu ci-dessus. D'o\`u la non nullit\'e de l'op\'erateur d'entrelacement normalis\'e comme en \cite{manuscripta}.
Pour les autres places,  la non nullit\'e vient de \cite{manuscripta} 3.5.
\section{Nombre de p\^oles}
Dans cette section on suppose qu'en toute place archim\'edienne le groupe est quasi d\'eploy\'e. On fixe une repr\'esentation cuspidale $\pi$ et on suppose qu'en toute place archim\'edienne la composante locale $\pi_{v}$ est dans le paquet de Langlands \`a l'int\'erieur du paquet d'Arthur; ce qui importe pour nous est qu'alors, par d\'efinition, $L(\rho_{v}\times\pi_{v},s)=L(\rho_{v}\times \pi^{\GL}_{v},s)$. Sous cette hypoth\`ese, la proposition de \ref{eisenstein} est vraie: l'holomorphie et la non nullit\'e des op\'erateurs d'entrelacement aux places archim\'ediennes r\'esulte de \cite{manuscripta} 3.5.

On \'ecrit $\pi^{\GL}\simeq \times_{(\rho',b')\in {\mathcal {S}}}{\Speh}(\rho',b')$, ce qui d\'efinit ${\mathcal S}$. On note ${\mathcal {S}}_{\rho}$ l'ensemble des entiers $b$ tel que $(\rho,b)\in {\mathcal {S}}$.
\begin{cor} Le nombre de p\^oles r\'eels positifs de la fonctions $L(\rho\times \pi,s)/L(\rho\times \pi,s+1)$ est inf\'erieur ou \'egal au cardinal du sous-ensemble de ${\mathcal {S}}_{\rho}$ form\'e des entiers $b$ tel que $b+2\notin {\mathcal {S}}_{\rho}$.
\end{cor}
Soit $s_{0}\in {\mathbb R}_{>0}$ un p\^ole de la fonciton $L(\rho\times \pi,s)/L(\rho\times \pi,s+1)$. D'apr\`es \ref{eisenstein} et ce qui pr\'ec\`ede l'\'enonc\'e, on sait que la s\'erie d'Eisenstein $E(\rho\times \pi,s)$ a aussi un p\^ole en $s=s_{0}$ et son r\'esidu est de carr\'e int\'egrable. Notons $\pi_{+}$ ce r\'esidu; on sait que c'est une repr\'esentation irr\'eductible (cf. \cite{image}) mais pour ce que l'on fait ici, il suffit d'en prendre une composante irr\'eductible. On \'ecrit le transfert tordu, $\pi_{+}^{\GL}$ de $\pi_{+}$ au bon groupe lin\'eaire: n\'ecessairement, $(\rho,2s_{0}-1)\in {\mathcal{S}}$ et 
$$
\pi_{+}^{\GL}\simeq {\Speh}(\rho,2s_{0}+1)\times _{(\rho',b')\in {\mathcal{S}}; (\rho',b')\neq (\rho,2s_{0}-1)}{\Speh}(\rho',b').
$$
Donc en particulier $(\rho,2s_{0}+1)\notin {\mathcal S}$. Cela prouve le corollaire.
\begin{rmq} Si on ne suppose pas que $\pi$ est cuspidal, alors le nombre de p\^oles r\'eels positifs de la fonction $L(\rho\times \pi,s)/L(\rho\times \pi,s+1)$ est inf\'erieur ou \'egal \`a $2\vert {\mathcal{S}}_{\rho}\vert$ cette borne pouvant tr\`es certainement \^etre atteinte.
\end{rmq}
En effet, il se peut tout \`a fait qu'en toute place $v$, le groupe soit quasi-d\'eploy\'e et la composante locale de $\pi$ soit dans le paquet de Langlands \`a l'int\'erieur du paquet d'Arthur. Si cette propri\'et\'e est v\'erifi\'ee, alors $L(\rho\times \pi,s)=L(\rho\times \pi^{\GL},s)$ et
$$
\frac{L(\rho\times \pi^{\GL},s)}{L(\rho\times \pi^{\GL},s+1)}=\prod_{(\rho',b')\in {\mathcal{S}}}\frac{L(\rho\times( \rho')^*,s-(b'-1)/2)}{L(\rho\times (\rho')^*,s+(b'+1)/2)}.
$$
On sait calculer les p\^oles de chaque membre du terme de droite; \'evidemment, ces p\^oles arrivent en des points demi-entiers et on ne peut exclure l'existence de $0$. On s'attend quand m\^eme qu'il existe une donn\'ees ${\mathcal{S}}$ tel que pour $(\rho',b'),(\rho'',b'')\in {\mathcal {S}}$, $L(\rho'\times (\rho'')^*,1/2)\neq 0$. Dans ce cas, la borne est de l'\'enonc\'e est bien atteinte.

\bibliography{fonctionsL}

\bibliographystyle{alpha}

\end{document}